\documentclass[12pt]{amsart}     

\usepackage{tkz-fct}

\usepackage[margin=1.15in]{geometry}

\usepackage[pagebackref=true,
]{hyperref} 
\hypersetup{colorlinks=true,linkcolor=blue,filecolor=magenta,citecolor=green,urlcolor=blue,final} 

 \usepackage{amssymb}

\usepackage{amsthm,amssymb,amsmath,amscd}
\usepackage{amsfonts,color}
\usepackage{latexsym}
\usepackage[all]{xy}
\usepackage{ulem}

\usepackage{mathptmx}

\newtheorem{thm}{Theorem}[section]
\newtheorem{lem}[thm]{Lemma}
\newtheorem{cor}[thm]{Corollary}
\newtheorem{prop}[thm]{Proposition}
\theoremstyle{definition}
\newtheorem{df}[thm]{Definition}
\theoremstyle{remark}
\newtheorem{rem}[thm]{Remark}

\newtheorem{example}[thm]{Example}

\numberwithin{equation}{section}

\newtheorem*{ack}{Acknowledgements}

\newcommand{\bC}{{\mathbb C}}

\newcommand{\bN}{{\mathbb N}}

\newcommand{\bR}{{\mathbb R}}
\newcommand{\bQ}{{\mathbb Q}}
\newcommand{\bZ}{{\mathbb Z}}
\newcommand{\bT}{{\mathbb T}}

\newcommand{\cE}{{\mathcal E}}
\newcommand{\cF}{{\mathcal F}}

\newcommand{\cO}{{\mathcal O}}

\newcommand{\Hom}{{\rm{Hom}}}

\newcommand{\wti}{\widetilde}

\newcommand{\lra}{\longrightarrow}

\newcommand{\Int}{{\rm Int}}
\newcommand{\Relint}{{\rm Relint}}

\newcommand{\ch}{{\rm ch}}

\newcommand{\td}{{\rm td}}
\newcommand{\mult}{{\rm mult}}

\newcommand{\bb}[1]{\mbox{$\mathbb{#1}$}}

\def\sig{\sigma}
\def\Sig{\Sigma}

\def\be{\begin{equation}}
\def\ee{\end{equation}}
\def\bt{\begin{thm}}
\def\et{\end{thm}}
\def\bc{\begin{cor}}
\def\ec{\end{cor}}
\def\br{\begin{rem}}
\def\er{\end{rem}}
\def\bp{\begin{prop}}
\def\ep{\end{prop}}
\def\bl{\begin{lem}}
\def\el{\end{lem}}
\def\bn{\begin{enumerate}}
\def\en{\end{enumerate}}
\def\bex{\begin{example}}
\def\eex{\end{example}}
\def\bd{\begin{df}}
\def\ed{\end{df}}

\usepackage{ textcomp }

\begin{document}                        


\title[Euler-Maclaurin formulae]{Equivariant toric geometry \\ and Euler-Maclaurin formulae \\ -- an overview --  }
\dedicatory{Dedicated to Lauren\c{t}iu P\u{a}unescu and Alexandru Suciu}

\author[S. E. Cappell ]{Sylvain E. Cappell}
\address{S. E. Cappell: Courant Institute, New York University, 251 Mercer Street, New York, NY 10012, USA.}
\email {cappell@cims.nyu.edu}

\author[L. Maxim ]{Lauren\c{t}iu Maxim}
\address{L. Maxim : Department of Mathematics, University of Wisconsin-Madison, 480 Lincoln Drive, Madison WI 53706-1388, USA, \newline
{\text and} \newline Institute of Mathematics of the Romanian Academy, P.O. Box 1-764, 70700 Bucharest, ROMANIA.}
\email {maxim@math.wisc.edu}

\author[J. Sch\"urmann ]{J\"org Sch\"urmann}
\address{J.  Sch\"urmann : Mathematische Institut,
          Universit\"at M\"unster,
          Einsteinstr. 62, 48149 M\"unster,
          Germany.}
\email {jschuerm@math.uni-muenster.de}

\author[J. L. Shaneson ]{Julius L. Shaneson}
\address{J. L. Shaneson: Department of Mathematics, University of Pennsylvania, 209 S 33rd St., Philadelphia, PA 19104, USA.}
\email {shaneson@sas.upenn.edu}



\date{\today}

\begin{abstract}
We survey recent developments in the study of torus equivariant motivic Chern and Hirzebruch characteristic classes of projective toric varieties, with applications to calculating equivariant Hirzebruch genera of torus-invariant Cartier divisors in terms of torus characters, as well as to general 
Euler-Maclaurin type formulae for full-dimensional simple lattice polytopes. We present recent results by the authors, emphasizing the main ideas and some key examples. This includes global formulae for equivariant Hirzebruch classes in the simplicial context proved by localization at the torus fixed points, a weighted versions of a classical formula of Brion, as well as of the Molien formula of Brion-Vergne.

Our Euler-Maclaurin type formulae provide generalizations to arbitrary coherent sheaf coefficients of the Euler-Maclaurin formulae of Cappell-Shaneson, Brion-Vergne, Guillemin, etc., via the equivariant Hirzebruch-Riemann-Roch formalism. Our approach, based on motivic characteristic classes, allows us, e.g., to obtain such Euler-Maclaurin formulae also for (the interior of) a face. We obtain such results also in the weighted context, and for Minkovski summands of the given full-dimensional lattice polytope. \end{abstract}

\maketitle





\section{Introduction. Historical overview}
Counting lattice points in full-dimensional polytopes is a problem with a long mathematical history. A natural generalization of this problem is to compute the sum of the values of a suitable continuous function $f$ at the lattice points contained in a polytope $P$. The resulting expression is called the {\it Euler-Maclaurin formula} for $f$ and $P$. 

\smallskip

As observed by Danilov \cite{D}, the lattice point counting problem is closely related to the \index{Riemann-Roch theorem} Riemann-Roch theorem \cite{BFM}, via the correspondence between lattice polytopes and toric varieties. An $n$-dimensional {toric variety} $X=X_\Sigma$ is an irreducible normal variety on which the complex affine $n$-torus $\bT\simeq (\bC^*)^n$ acts with an open orbit, e.g., see \cite{CLS, D,F1}. Toric varieties arise from combinatorial objects $\Sigma \subset N \otimes \bR \simeq \bR^n$ called fans, which are collections of cones in a lattice $N \simeq \bZ^n$. Here $N$ corresponds to one-parameter subgroups of $\bT$. 
Let $M \simeq \bZ^n$ be the character lattice of $\bT$. 
From a full-dimensional lattice polytope $P\subset M \otimes \bR\simeq \bR^n$ one constructs (via the associated \index{inner normal fan} inner normal fan $\Sigma_P$ of $P$) a {toric variety} $X_P:=X_{\Sig_P}$, together with an ample Cartier divisor $D_P$, so that  the number of lattice points of $P$, i.e., points of $M \cap P$, is computed by the holomorphic Euler characteristic $\chi(X_P,\cO(D_P))$.  The Riemann-Roch theorem expresses the latter in terms of  the Chern character of $\cO(D_P)$ and the Baum-Fulton-MacPherson homology \index{Todd class} Todd class $\td_*(X_P)$ of the toric variety $X_P$, thus reducing  the lattice point counting problem to a characteristic class computation. 


In \cite{MS1, MS2}, the second and third author computed the \index{motivic Chern class} {\it motivic Chern class} $mC_y$ and, resp., homology  {\it Hirzebruch classes} $T_{y*}$ \cite{BSY}  of (possibly singular) toric varieties.  
In particular, they obtained new, or recovered well-known, formulae for the Baum-Fulton-MacPherson Todd  classes $\td_*(X)=T_{0*}(X)$.
By taking advantage of the torus-orbit decomposition and  the motivic properties of the homology  
Hirzebruch classes, one can express the latter in terms of the (dual) 
 Todd classes of closures of orbits.  As a consequence, by generalizing Danilov's observation, one shows that the weighted lattice point counting, where each point in a face $E$ of the polytope 
$P$ carries the weight $(1+y)^{\dim(E)}$, amounts to the computation of the Hirzebruch class $T_{y*}(X_P)$ of the associated toric variety. 
\smallskip

An Euler-Maclaurin formula  relates the sum $\sum_{m\in P \cap M} f(m)$ of the values of a suitable function $f$ at the lattice points in a lattice polytope $P \subset M_\bR:=M\otimes\bb{R}$ to integrals over the polytope and/or its faces. In \cite{CMSSEM}, we consider $f$ to be a polynomial on $M_\bR$, or an exponential function $f(m)=e^{\langle m, z \rangle}$, or products of these two types, where $\langle \cdot , \cdot \rangle: M\times N\to \bb{Z}$ is the canonical pairing and $z \in N_\bC:=N \otimes_\bZ \bC=\Hom_\bR(M_\bR, \bC)$. 

Khovanskii and Pukhlikov \cite{KP} obtained an
Euler-Maclaurin formula for the sum of the values of a polynomial over the lattice points in {regular} lattice polytopes (corresponding to smooth projective toric varieties). The first substantial advance for non-regular polytopes was a different type of Euler-Maclaurin formula for {simple} polytopes achieved by the first and fourth authors in \cite{CS2,S}. 
A few years later, the Khovanskii-Pukhlikov formula was extended to {simple} lattice polytopes by Brion-Vergne \cite{BrV1, BrV2}, the latter using the {\it equivariant Hirzebruch-Riemann-Roch} theorem for the corresponding complete simplicial toric varieties, together with localization techniques in equivariant cohomology. Other Euler-Maclaurin type formulae 
were obtained by Guillemin \cite{G} by using methods from symplectic geometry and geometric quantization, by Karshon et al. \cite{KSW} by combinatorial means, etc. 
While most of above-mentioned Euler-Maclaurin formulae are obtained by integrating the function $f$ over a dilation of the polytope, the Cappell-Shaneson approach involves a summation over the faces of the polytope of integrals (over such faces) of linear differential operators with constant coefficients, applied to the function; see also \cite{BrV3} for some formulae of this type. Later on, these Euler-Maclaurin formulae have been extended (even for non-simplicial toric varieties, resp.,  non-simple lattice polytopes)
by Berline--Vergne \cite{BeV1, BeV2} and Garoufalidis--Pommersheim \cite{GP}, together with Fischer--Pommersheim \cite{FP}, to local formulae satisfying a {\it Danilov condition} (see also the work of Pommershein--Thomas \cite{PT} for such formulae for the Todd class in the non-equivariant context). In the geometric context of toric varieties, all of these Danilov type formulae are tied to the birational invariance of the (equivariant) Todd class of the structure sheaf. However, this birational invariance property does not apply to the more general context considered in our paper \cite{CMSSEM}.

\medskip

In this survey, we overview some of the recent results obtained in \cite{CMSSEM}, where  we consider $\bT$-equivariant versions $mC_y^\bT$ and $T_{y*}^\bT:=\td_*^\bT \circ mC_y^\bT$
of the motivic Chern and, resp.,  Hirzebruch characteristic classes of \cite{BSY}, and extend the formulae from \cite{MS1, MS2} to the equivariant setting. 
 As in \cite{CMSSEM}, we also elaborate here on the relation (cf. also \cite{BrV2}) between the  equivariant toric  geometry via the (equivariant) Hirzebruch-Riemann-Roch (abbreviated HRR for short) and Euler-Maclaurin type formulae for simple lattice polytopes (corresponding to simplicial toric varieties).
 
Our main results from \cite{CMSSEM} provide generalizations to {\it arbitrary coherent sheaf coefficients}, which for natural choices related to the toric variety (or the polytope) give uniform geometric proofs of the Euler-Maclaurin formulae of Brion-Vergne and, resp., Cappell-Shaneson, via the equivariant Hirzebruch-Riemann-Roch formalism. Our approach, based on motivic characteristic classes, allows us to obtain such Euler-Maclaurin formulae also for (the interior of) a face, as well as for the polytope with several facets (i.e., codimension one faces) removed, e.g., for the interior of the polytope. 
Moreover, we prove such results also in the weighted context, as well as for $\bN$-Minkowski summands of the given full-dimensional lattice polytope (corresponding to globally generated torus invariant Cartier divisors in the toric context). Similarly, some of these results are extended to local Euler-Maclaurin formulae for the tangent cones at the vertices of the given  full-dimensional lattice polytope (fitting with localization at the torus fixed points in equivariant $K$-theory and equivariant (co)homology).


In what follows we give a brief description of some of our main results from \cite{CMSSEM}, emphasizing the main ideas and some key examples. We invite the interested reader to consult \cite{CMSSEM} for more details and complete proofs of these results. For simplicity, in this survey we formulate the $K$-theoretical results for a projective toric variety, whereas for the cohomological counterparts we further restrict to the simplicial context (e.g., the toric variety associated to a simple full-dimensional lattice polytope); see \cite{CMSSEM} for more general statements.

\section{Equivariant Hirzebruch-Riemann-Roch}

We first introduce some notations. 

\subsection{Rational equivariant cohomology}
Let $X=X_\Sig$ be an $n$-dimensional projective simplicial  toric variety with fan $\Sig \subset N_\bR=N \otimes \bR$ and torus $\bT=T_N$.
Denote by $H^*_\bT(X;\bb{Q})$ the (Borel-type) rational \index{equivariant cohomology} equivariant cohomology of $X$, and note that for a point space one has  
\[ H^*_\bT(pt;\bQ)\simeq \bQ[t_1,\ldots,t_n]=:(\Lambda_\bT)_{\bQ}.\]
Let $M$ be the dual lattice of $N\simeq \bZ^n$.
Viewing characters $m\in M$   (resp., $\chi^m\in \bZ[M]\simeq K_0^\bT(pt)$) of $\bT$ as $\bT$-equivariant line bundles $\bC_{\chi^m}$ over a point space $pt$ gives an isomorphism $M\simeq Pic_\bT(pt)$. Taking the first equivariant Chern class 
$c^1_\bT$ (or the dual $-c^1_\bT$) gives an  isomorphism
\[
c=c^1_\bT,  \:\text{resp.}, \:  s=-c^1_\bT :\:
M \simeq H_\bT^2(pt;\bb{Z}).
\]
Hence, upon choosing a basis $m_i$ ($i=1,\dots,n$) of $M\simeq \bb{Z}^n$, one has that $H_\bT^*(pt;\bb{Q})=(\Lambda_\bT)_{\bb{Q}}\simeq \bb{Q}[t_1,\dots,t_n],$  with
$t_i=\pm c^1_\bT(\bC_{\chi^{m_i}})$ for $i=1,\dots,n$. Moreover, $H^*_\bT(X;\bb{Q})$ can be described as a  $H_\bT^*(pt;\bb{Q})=(\Lambda_\bT)_{\bb{Q}}$-algebra. This fact plays an important role for proving Euler-Maclaurin type formulae. In fact, we will be working with the completions $\widehat{H}_\bT^*(X;\bb{Q}) := \prod_{i\geq 0} \: H^i_\bT(X;\bb{Q})$ and $(\widehat{\Lambda}_\bT)_{\bb{Q}}\simeq \bQ[[t_1,\ldots,t_n]]$ of these rings. It is important to note that the equivariant Chern character $\ch^\bT$ and equivariant Todd homology class transformation $\td_*^\bT$ of Edidin-Graham \cite{EG} and Brylinski-Zhang \cite{BZ} take values in an {\it analytic subring} (cf. \cite[Proposition 5.17]{CMSSEM})
$$(H^*_\bT(X;\bb{Q}))^{an}\subset  \widehat{H}^*_\bT(X;\bb{Q}) \:,$$
with $ \bb{Q}\{t_1,\dots,t_n\}\simeq (H^*_\bT(pt;\bb{Q}))^{an}=:(\Lambda^{an}_\bT)_{\bb{Q}}\subset (\widehat{\Lambda}_\bT)_{\bb{Q}}$ the subring of {\it convergent power series} (around zero) with rational coefficients, i.e., after pairing with $z\in N_\bC$ one gets a convergent power series {\it function} in $z$ around zero, whose corresponding Taylor polynomials have rational coefficients. Here we use equivariant Poincar\'{e} duality $\widehat{H}_*^\bT(X;\bb{Q}) \simeq \widehat{H}^*_\bT(X;\bb{Q})$ between equivariant homology and cohomology for simplicial toric varieties, to view these classes in equivariant cohomology.

\subsection{Generalized equivariant Hirzebruch-Riemann-Roch formula and applications}
Let $X=X_\Sig$ be a projective toric variety, with a $\bT$-equivariant coherent sheaf $\cF$. In this section, we do not need to assume that $X$ is simplicial. The cohomology spaces $H^i(X;\cF)$ are finite dimensional $\bT$-representations, vanishing for $i$ large enough. Using the corresponding $\bT$-eigenspaces $H^i(X;\cF)_{\chi^m}$ as in  \cite[Proposition 1.1.2]{CLS}, on which $t \in \bT$ acts as multiplication by $\chi^m(t)$, the (cohomological)  {Euler characteristic} of $\cF$ is defined by
\be\label{ie1}
\chi^\bT(X,\cF)=\sum_{m \in M} \sum_{i=0}^n (-1)^i \dim_\bC H^i(X;\cF)_{\chi^m} \cdot e^{c(m)} \in (\Lambda^{an}_\bT)_{\bb{Q}}
\subset  (\widehat{\Lambda}_\bT)_{\bb{Q}}\:.
\ee
For $\cE$, resp., $\cF$, a $\bT$-equivariant vector bundle, resp., coherent sheaf on $X$, one then has the following \index{equivariant Hirzebruch-Riemann-Roch formula} {\it equivariant Hirzebruch-Riemann-Roch formula} (see \cite[Corollary 3.1]{EG}): 
\be\label{eHRRi}
\chi^\bT(X,\cE \otimes \cF)=\int_X \ch^\bT(\cE) \cap \td^\bT_*([\cF]),
\ee
where $\int_X:\widehat{H}_*^\bT(X;\bb{Q}) \to \widehat{H}_*^\bT(pt;\bb{Q})=(\widehat\Lambda_\bT)_\bQ$ is the equivariant pushforward for the constant map $X \to pt$. 

Formula \eqref{eHRRi} is extended in \cite{CMSSEM} to a \index{generalized equivariant Hirzebruch-Riemann-Roch formula} {\it generalized equivariant Hirzebruch-Riemann-Roch formula}, cf. \cite[Theorem 3.14]{CMSSEM} (where it is proved more generally for closed algebraic subsets of $X$ defined by $\bT$-invariant closed subsets):
\be\label{gHRRi}
\begin{split}
\chi^\bT_y(X,\cO_{X}(D)):=
\sum _{p=0}^{n} \chi^\bT(X,\widehat{\Omega}^{p}_X \otimes  \cO_{X}(D)) \cdot y^p
=\int_{X} \ch^\bT(\cO_{X}(D)) \cap T^\bT_{y*}(X),
\end{split}
\ee
with $D$ a $\bT$-invariant Cartier divisor on $X$, $\ch^\bT$ the equivariant Chern character,  and $\widehat{\Omega}^{p}_X$ the sheaf of Zariski 
$p$-forms on $X$. Here  we use the following explicit description of the equivariant motivic Chern class $mC_y^\bT(X)$ and, resp.,  equivariant  Hirzebruch class $T^\bT_{y*}(X)$ of $X$, obtained in \cite[Proposition 3.5]{CMSSEM}):
\be\label{4} mC^\bT_y(X)=\sum_{p=0}^{\dim(X)} [\widehat{\Omega}^p_X]_\bT \cdot y^p \in K_0^\bT(X)[y]
\quad \text{and} \quad 
{T}^\bT_{y*}(X)=\sum_{p=0}^{\dim(X)} \td_*^\bT([\widehat{\Omega}^p_X]_\bT) \cdot y^p.
\ee
The reader unfamiliar with these characteristic class notions, can take the above formulae from \eqref{4} as their definitions.

Let now $P$ be a full-dimensional \index{simple lattice polytope} lattice polytope in $M_\bR\simeq \bR^n$ with associated toric variety $X=X_P$ with torus $\bT$, inner normal fan $\Sigma=\Sig_P$ and ample Cartier divisor $D=D_P$. 
As a consequence of \eqref{gHRRi}, we obtain the following weighted formula (see \cite[Corollary 3.15]{CMSSEM}):
\be\label{ic1}
\chi^\bT_y(X,\cO_{X}(D))= \sum_{E \preceq P} (1+y)^{\dim(E)} \cdot \sum_{m \in \Relint(E) \cap M} e^{s(m)},
\ee
where the first sum is over the faces $E$ of $P$ and $\Relint(E)$ denotes the relative interior of the face $E$. Let us also explain the use of $s(m)$ instead of $c(m)$ in the above formula, in the simple case when $y=0$.  
Formula \eqref{ic1} is then based on the vanishing of higher cohomology of $\cO_{X}(D)$, together with (e.g., see  \cite[Prop.4.3.3]{CLS})
\be\label{am1}
\Gamma(X;\cO_{X}(D)) = \bigoplus_{m \in P \cap M} \bC \cdot \chi^m \subset \bC[M]=\Gamma(\bT, \cO_\bT),
\ee
with $\bT$ acting on $\Gamma(\bT, \cO_\bT)$ as follows: if $t \in \bT$ and $f \in \Gamma(\bT, \cO_\bT)$, then $t \cdot f \in \Gamma(\bT, \cO_\bT)$ is given by $p \mapsto f(t^{-1} \cdot p)$, for $p \in \bT$ (see \cite[pag.18]{CLS}), so that $\Gamma(X;\cO_{X}(D))_{\chi^{-m}} = \bC \cdot \chi^{m}$ for all $m\in P \cap M$.


Let us next consider a globally generated $\bT$-invariant Cartier divisor $D'$ on $X=X_\Sig$, with associated (not necessarily full-dimensional) lattice polytope $P_{D'} \subset M_\bR$. Let $X_{D'}$ be the toric variety of the lattice polytope $P_{D'}$, defined via the corresponding {generalized fan} $\Sig'$ as in \cite[Prop.6.2.3]{CLS}. There is a proper toric morphism $f:X \to X_{D'}$, induced by the corresponding lattice projection $N \to  N_{D'}$ given by dividing out by the minimal cone of the generalized fan of $P_{D'}$. In particular, $f\colon X \to X_{D'}$ is a \index{toric fibration} toric fibration. For $\sig'$ a cone in the generalized fan $\Sig'$  of $P_{D'}$, let
\be\label{bl1}
d_\ell(X/\sig'):=\vert \Sigma_\ell(X/\sig') \vert,\ee with  
\be\label{bl2} \Sigma_\ell(X/\sig'):=\{\sig \in \Sig \mid  f(O_\sig)=O_{\sig'}, \ \ell=\dim(O_\sig) - \dim(O_{\sig'})\},\ee
where $O_\sig$ ($\sig \in \Sig$) and $O_{\sig'}$ ($\sig' \in \Sig'$) are $\bT$-, and resp. $\bT'$-orbits,
and $\vert - \vert$ denotes the cardinality of a finite set.
If $E$ is the face of $P_{D'}$ corresponding to $\sig' \in \Sig'$, we denote these multiplicities by $d_\ell(X/E)$.
Then we have the following generalization of formula \eqref{ic1} (see \cite[Corollary 3.17]{CMSSEM} for a more general statement):
\be\label{ic2}
\chi^\bT_y(X,\cO_{X}(D'))= \sum_{E \preceq P_{D'}} \left( \sum_{\ell \geq 0} 
(-1)^\ell \cdot d_\ell(X/E) \cdot (1+y)^{\ell + \dim(E)}\right)  \cdot  \sum_{m \in \Relint(E) \cap M} e^{s(m)}.
\ee
By forgetting the $\bT$-action (i.e., setting $s(m)=0$ for all $m \in M$), 
 we get a weighted lattice point counting for lattice polytopes associated to globally generated $\bT$-invariant Cartier divisors:
 \be\label{fort}
 \chi_y(X,\cO_{X}(D'))= \sum_{E \preceq P_{D'}} \left( \sum_{\ell \geq 0} 
(-1)^\ell \cdot d_\ell(X/E) \cdot (1+y)^{\ell +\dim(E)}  \right) \cdot  \vert \Relint(E) \cap M \vert .
\ee

\section{Localized equivariant motivic Chern and Hirzebruch classes}
In \cite[Section 4]{CMSSEM}, we apply localization techniques in $\bT$-equivariant $K$-theory for toric varieties (due to Brion-Vergne \cite{BrV2}), and $\bT$-equivariant cohomology for simplicial toric varieties (due to Brylinski-Zhang \cite{BZ} in the more general equivariant homology context), for the calculation of the $\bT$-equivariant motivic Chern (and Hirzebruch classes) of (simplicial) projective toric varieties.

Let $X=X_\Sig$ be an $n$-dimensional projective toric variety with torus $\bT=T_N$, so the fixed-point set $X^\bT \neq \emptyset$. Let $x_\sig \in X^\bT$ be the fixed point for $\sig \in \Sig(n)$, with $U_\sig$ the corresponding $\bT$-invariant open affine variety. Consider the multiplicative subset $S\subset \bZ[M]=K_0^\bT(pt)$ generated by the elements $1-\chi^m$, for $0\neq m \in M$. Then we have the localization isomorphism of Brion-Vergne, $K_0^\bT(X^\bT)_S \simeq K_0^\bT(X)_S$, induced from the inclusion of fixed points.
The projection map 
$$pr_{x_\sig} : K_0^\bT(X)_S \simeq K_0^\bT(X^\bT)_S \to  K_0^\bT(x_\sig)_S=\bZ[M]_S,$$ 
can be calculated, after restriction to $U_\sig$, as $pr_{x_\sig}=\mathbb{S} \circ \chi^\bT_\sig$, with $$\chi^\bT_\sig:K_0^\bT(U_\sig) \to \bZ[M]_{\rm sum} \subset \bZ[[M]]$$  the local counterpart of the equivariant Euler characteristic \eqref{ie1} (using $\chi^m\in \bZ[M]$ instead of $e^{c(m)}$),  and $\mathbb{S}$ the corresponding summation map as introduced by Brion-Vergne \cite{BrV2}. 
We then have (cf. \cite[formula (155)]{CMSSEM} for a more general version):
\be\label{imclocal}
\chi_\sig^\bT(mC^\bT_y(X)\vert_{U_\sig})=
\sum_{\tau \preceq \sigma} (1+y)^{\dim(O_\tau)} 
\sum_{m \in {\rm Relint}(\sigma^\vee \cap \tau^\perp) \cap M} \chi^{-m} \in  \bb{Z}[M]_{\rm sum} \otimes_\bZ \bZ[y],\ee
where the first sum is over the faces of $\sig$, with $\sigma^\vee$ is the dual cone and $\tau^\perp \subset M_\bR$ is the orthogonal of $\tau \subset N_\bR$ with respect to the canonical pairing. 
For $y=0$, this specializes to $\chi_\sig^\bT(\cO_X\vert_{U_\sig})$, since $X$ has rational singularities, so that $mC_0^\bT(X)=[\cO_X]_\bT\in K_0^\bT(X)$. As a consequence, we get the following weighted version of Brion's formula (cf. \cite[Corollary 4.8]{CMSSEM}).
\bc\label{iwBr}
Let $P$ be a full-dimensional lattice polytope with associated projective toric variety $X=X_P$ and ample Cartier divisor $D=D_P$. For each vertex $v$ of $P$, consider the cone $C_v={\rm Cone}(P \cap M -v)=\sigma_v^\vee$, with faces $E_v={\rm Cone}(E \cap M -v)$ for $v\in E$. 
Then the following identity holds in $\bZ[M]_S \otimes_\bZ \bZ[y]$:
\be\label{if94} \chi^\bT(X, mC^\bT_y(X) \otimes \cO_X(D))= \sum_{v \ \text{\rm vertex}} \chi^{-v} \cdot \mathbb{S} \left( \sum_{v \in E \preceq P} (1+y)^{\dim(E)} \cdot \sum_{m \in \Relint(E_v) \cap M} \chi^{-m} \right).
\ee
\ec
Brion's formula \cite{Br1} is obtained from \eqref{if94} by specializing to $y=0$.

In the case of a {\it smooth} cone $\sigma$ with $m_{\sig,i}$, $i=1,\ldots,n$, the minimal generators of $\sig^\vee$, formula \eqref{imclocal} becomes
$$
\chi_\sig^\bT(mC^\bT_y(X)\vert_{U_\sig})=\prod_{i=1}^n \left(1+(1+y) \cdot \sum_{k\geq 1} 
(\chi^{-m_{\sig,i}})^k \right),
$$
hence  
\be\label{sumsmmc}
\begin{split}
\mathbb{S}(\chi_\sig^\bT(mC^\bT_y(X)\vert_{U_\sig}))=\prod_{i=1}^n \left( 1+(1+y) \cdot \frac{\chi^{-m_{\sig,i}}}{1-\chi^{-m_{\sig,i}}}  \right) 
=\prod_{i=1}^n \frac{1+y \cdot \chi^{-m_{\sig,i}}}{1-\chi^{-m_{\sig,i}}}
\in \bZ[M]_S.
\end{split}
\ee

In the case of a \index{simplicial cone} {\it simplicial} cone, we get similar explicit formulae by using a Lefschetz type variant 
$tr^{\bT'}_\sig$ of the Euler characteristic $\chi_\sig^\bT$, and a corresponding summation map (see \cite[formula (144)]{CMSSEM} for more details). This Lefschetz type variant is an adaptation of the classical linear algebra formula (with trace denoted by $tr$)
$$\dim V^G=\frac{1}{\vert G \vert} \sum_{g \in G} tr(g:V \to V),$$
for a finite linear group action of $G$ on a finite dimensional complex vector space $V$.
 Let $\sig \in \Sig(n)$ be an $n$-dimension simplicial cone with $u_1,\ldots, u_n \in N=N_\sig$ the generators of the the rays $\rho_j \in \sig(1)$ of the cone $\sigma$, $j=1,\ldots,n$. Let $N'=N'_\sig$ be the finite index sublattice of $N$ generated by $u_1,\ldots, u_n$, and consider $\sig \in N'_{\bR}=N_\bR$ so that it is smooth with respect to the lattice $N'$. With $\bT$, $\bT'$ the corresponding $n$-dimensional tori of the lattices $N$, resp., $N'$, the inclusion $N' \hookrightarrow N$ induces a toric morphism $\pi: U'_\sig \to U_\sig$ of the associated affine toric varieties. Let $G_\sig$ be the finite kernel of the epimorphism $\pi:\bT'\to \bT$, so that $U'_\sig/G_\sig \simeq U_\sig$. 
Let $m'_{\sig,1},\ldots, m'_{\sig,n}$ be the dual basis in the dual lattice $M'=M_{\sig}'$ of $N'$, with corresponding characters $a_{\rho_j}:G_\sig \subset \bT'\simeq (\bC^*)^n \to \bC^*$ of $G_\sig$ given by the projection onto the $\rho_j$-th factor.    With these notations, we have (see \cite[formula (152)]{CMSSEM}):
\be\label{isumsmusd}
\mathbb{S}(\chi_\sig^\bT(mC^\bT_y(X)\vert_{U_\sig}))=\frac{1}{\vert G_\sig \vert} \sum_{g \in G_\sig}  \prod_{i=1}^n \frac{1+y \cdot a_{\rho_i}(g^{-1}) \cdot \chi^{-m'_{\sig,i}}}{1-a_{\rho_i}(g^{-1}) \cdot \chi^{-m'_{\sig,i}}}.
\ee
Here, $a_{\rho_i}(g^{-1}) \in \bC^*$ are the traces  of the action of $g^{-1}$ on the $1$-dimensional $\bT'$-representations corresponding to the characters $\chi^{-m'_{\sig,i}}$. See \cite[Example 4.4]{CMSSEM} for more details.
For $y=0$, \eqref{isumsmusd} specializes to the \index{Molien formula} {\it Molien formula} of Brion-Vergne \cite{BrV2}.

\medskip

We next assume that $X$ is a projective simplicial toric variety, and discuss some cohomological counterparts of the above localization formulae. 
Let $L\subset (\Lambda_\bT)_\bQ=H^{*}_\bT(pt;\bQ)$ be the multiplicative subset  generated by the elements $\pm c(m)$, for $0\neq m \in M$. With $x_\sig \in X^\bT$ a fixed point corresponding to $\sig \in \Sig(n)$, there is an associated cohomological localization map  at $x_\sig $, 
\be\label{locmap}\frac{ i_\sig^*}{Eu^\bT_X(x_\sig)}:\widehat{H}^{*}_\bT(X;\bQ)_L 
\lra \widehat{H}^{*}_\bT(x_\sig;\bQ)_L = L^{-1}(\widehat{\Lambda}_\bT)_\bQ,\ee
with $i_\sig:\{x_\sig\} \hookrightarrow X$ the inclusion map, and $Eu^\bT_X(x_\sig)$
 the generalized Euler class  of the fixed point $x_\sig$ in $X$ defined by $$0\neq Eu^\bT_X(x_\sig):=i_\sig^* \left( \mult(\sig) \cdot \prod_{\rho \in \sig(1)} [D_\rho]_\bT\right) \in \bQ\cdot L.$$ Here, $[D_\rho]_\bT$ denotes the equivariant fundamental class of the $\bT$-invariant divisor $D_\rho$ corresponding to the ray $\rho \in \Sig(1)$, and $\mult(\sig)=|G_\sig|$ is the multiplicity of $\sig$. 
 Moreover, if  
 $\int_X:\widehat{H}^*_\bT(X;\bb{Q})_L \to \widehat{H}^*_\bT(pt;\bb{Q})_L=L^{-1}(\widehat\Lambda_\bT)_\bQ$ is the equivariant Gysin map (or, equivalently, the equivariant pushforward) for the constant map $X \to pt$, then 
\be\label{f108n}
\int_X = \sum_{\sig \in \Sig(n)} \frac{i_\sig^*}{Eu^\bT_X(x_\sig)}: \widehat{H}^{*}_\bT(X;\bQ)_L  \to \widehat{H}^{*}_\bT(pt;\bQ)_L.
\ee
 See \cite[Proposition 4.17]{CMSSEM} for more details.

These $K$-theoretic and cohomological localization maps are compatible with the equivariant Todd class transformation of Edidin-Graham \cite{EG0} (and Brylinski-Zhang \cite{BZ}), in the following sense:
\bp\label{ipr47}
Let $\cF$ be a $\bT$-equivariant coherent sheaf on the projective simplicial toric variety $X=X_\Sig$, and let $x_\sig \in X^\bT$ be a given fixed point of the $\bT$-action. Then:
 \be\label{if98}
\td_*^\bT([\cF])_{x_\sig}:= \frac{ i_\sig^*\td_*^\bT([\cF])}{Eu^\bT_X(x_\sig)}  = \ch^\bT (( \mathbb{S} \circ \chi_\sig^\bT)(\cF))
\in  L^{-1}(\Lambda^{an}_\bT)_\bQ \subset L^{-1}(\widehat{\Lambda}_\bT)_\bQ,
 \ee
with $\ch^\bT:\bZ[M]_S \to L^{-1}(\Lambda^{an}_\bT)_\bQ$ induced by the $\bT$-equivariant Chern character on a point space. Moreover, the cohomological Euler characteristic of $\cF$ can be calculated via localization at the $\bT$-fixed points as:
\be\label{f104}
\chi^\bT(X,\cF) 
=\sum_{\sig \in \Sigma(n)} \td_*^\bT([\cF])_{x_\sig} \in L^{-1}(\Lambda^{an}_\bT)_\bQ \subset L^{-1}(\widehat{\Lambda}_\bT)_\bQ \: .
\ee
\ep

For the simplicial cone $\sig \in \Sig(n)$ corresponding to $x_\sig \in X^\bT$, formula \eqref{isumsmusd} implies by applying the Chern character the following:
\be\label{isum}
T^\bT_{y*}(X)_{x_\sig}:=\td_*^\bT([mC_y^\bT(X)])_{x_\sig} =\frac{1}{\vert G_\sig \vert} \sum_{g \in G_\sig}  \prod_{i=1}^n \frac{1+y \cdot a_{\rho_i}(g^{-1}) \cdot e^{-c(m'_{\sig,i})}}{1-a_{\rho_i}(g^{-1}) \cdot e^{-c(m'_{\sig,i})}},
\ee
By specializing \eqref{isum} to $y=0$, we get a formula for the localized Todd class  $\td_*^\bT(X)_{x_\sig}= T^\bT_{0*}(X)_{x_\sig}$, similar to \cite{BrV2}.

\section{Equivariant Hirzebruch classes}
The characteristic class formulae of \cite{MS1, MS2} for the motivic Chern and Hirzebruch classes are extended to the equivariant setting in \cite[Section 3]{CMSSEM} by using the global Cox construction \cite{Cox} which we now recall. 

\subsection{Cox construction}
For each ray $\rho \in \Sig(1)$ in the fan $\Sig$ of the projective simplicial toric variety $X=X_\Sigma$, denote by $u_{\rho}$ the corresponding ray generator. Let $r=| \Sig(1) |$ be the number of rays in the fan $\Sig$.
Using the fact that $N \simeq \Hom_{\bZ}(\bC^*,\bT)$ is identified with the one-parameter subgroups of $\bT$, define the map of tori $$\gamma:\widetilde{\bT}:=(\bC^*)^r \lra \bT \quad \text{by} \quad (t_{\rho})_{\rho} \mapsto \prod_{\rho \in \Sig(1)} u_{\rho}(t_{\rho}),$$ 
and let $G:=\ker (\gamma)$. Let $Z(\Sig) \subset \bC^r$ be the variety defined by the monoidal ideal generated by the elements $\hat{x}_{\sig}:=\prod_{\rho \notin \sig(1)} x_{\rho}$,  for $\sig \in \Sig$,
with $(x_{\rho})_{\rho \in \Sig(1)}$ the coordinates on $\bC^r$. Then the variety 
$W:=\bC^r \setminus Z(\Sig)$
is a toric manifold, and there is a toric morphism $\pi:W \to X.$
The group $G$ acts on $W$ by the restriction of the diagonal action of $(\bC^*)^{r}$, and the toric morphism $\pi$ is constant on $G$-orbits. Moreover, Cox \cite{Cox}  proved that if $X=X_{\Sig}$ is a projective simplicial toric variety,
then $X$ is the  geometric quotient $X=W/G$.

Let $$a_{\rho}:(\bC^*)^{r} \to \bC^*$$ be the projection onto the $\rho$-th factor. 
For a cone $\sigma \in \Sigma$, let $G_{\sig}$ be defined by
\be\label{stab}
\begin{split}
G_{\sig}:=\{g \in G \ | \ a_{\rho}(g)=1, \forall \rho \notin \sig(1) \}
\simeq \{(t_{\rho})_{\rho \in \sig(1)}  \ | \ t_{\rho}\in \bC^*, \prod_{\rho \in \sig(1)} u_{\rho}(t_{\rho})=1 \},
\end{split}
\ee
so $G_{\sig}$ 
depends only on $\sig$. Moreover, for $\sigma \in \Sigma(n)$ a top-dimensional cone, the last description of $G_\sig$ coincides with the definition of $G_\sig$ used in the local quotient description needed for formula \eqref{isumsmusd}. Note also that the corresponding characters $a_\rho$ in the global and local situations are identified. If $\tau \preceq \sig$ is a face of $\sig$, then by \eqref{stab}, one gets
$$\tau \preceq \sig  \Longrightarrow G_{\tau} \subseteq G_{\sig} \subset G.$$
Consider next the finite set $$G_{\Sig}:=\bigcup_{\sig \in \Sig} G_{\sig}=\bigcup_{\sig \in \Sig(n)} G_{\sig}.$$ 

\subsection{Equivariant Hirzebruch classes  of simplicial toric varieties}
The equivariant Hirzebruch class ${T}^\bT_{y*}(X)$ of a projective simplicial toric variety is computed by the following (cf. \cite[Theorem 3.22]{CMSSEM}):
\bt\label{thmhi} Let $X:=X_\Sig$ be an $n$-dimensional simplicial projective toric variety with fan $\Sig$. Then 
\be\label{ieHirz0}
{T}^\bT_{y*}(X)= (1+y)^{n-r} \cdot \sum_{g \in G_{\Sig}}  \prod_{\rho \in \Sig(1)} \frac{ F_{\rho} \cdot 
\big( 1+y  \cdot a_{\rho}(g)  \cdot e^{-F_{\rho}}\big)}{1-a_{\rho}(g) \cdot e^{-F_{\rho}}}  \in \widehat{H}^*_\bT(X;\bQ)[y] \:,
\ee
with $F_\rho=[D_\rho]_\bT$ denoting the equivariant fundamental class of the $\bT$-invariant divisor $D_\rho$ corresponding to the ray $\rho \in \Sig(1)$.
\et

For $y=0$, with $T_{0*}^\bT(X)=\td_*^\bT(X)$ the \index{equivariant Todd class} equivariant Todd class of $X$, formula \eqref{ieHirz0} specializes to the classical counterpart of Brion-Vergne \cite{BrV2} for the equivariant Todd class of $X$. A more general statement is obtained in \cite[Theorem 3.28]{CMSSEM}, for the equivariant Hirzebruch classes of complements of $\bT$-invariant divisors in $X$.

Theorem \ref{thmhi} is proved in \cite[Subsection 3.3]{CMSSEM} by using 
the equivariant Lefschetz-Riemann-Roch theorem of Edidin-Graham \cite{EG} for the geometric quotient $\pi:W \to X=W/G$. This can be seen as the global version of the previously discussed local Lefschetz-type arguments. These two pictures fit geometrically via the following identifications. 
For $\sig \in \Sigma(n)$ a top-dimensional cone, let $U_\sig$ be the $\bT$-invariant open affine subset of $X$ containing the corresponding $\bT$-fixed point $x_\sig$. We have $\pi^{-1}(U_\sig)\simeq \bC^n \times (\bC^*)^{r-n}$, with $\wti{\bT}\simeq \bT' \times (\bC^*)^{r-n}$ acting on the respective factors, the factor $\bC^n$ corresponding to the rays of $\sigma$, and $(\bC^*)^{r-n}$ acting freely by multiplication on itself. Similarly, $G\simeq G_\sigma \times (\bC^*)^{r-n}$, with $G_\sig \subset \bT'$ the finite subgroup introduced before. So, above $U_\sig$, $\pi$ can be factorized as a composition of the free quotient $\wti{\pi}:\bC^n \times (\bC^*)^{r-n} \to \bC^n$ by the $(\bC^*)^{r-n} $-action, followed by a finite quotient map $\pi=\pi_\sig:\bC^n \to \bC^n/{G_\sig}=U_\sig$.

Let us now indicate how localization can be used for a second proof of Theorem \ref{thmhi} (see \cite[Subsection 4.3.2]{CMSSEM} for a more general result).
\begin{proof}[Proof of Theorem \ref{thmhi}]
One starts by noticing that, for a projective simplicial toric variety, the map  
$$\bigoplus_{\sig \in \Sig(n)} i_\sig^* : \widehat{H}^*_\bT(X;\bQ) \hookrightarrow \bigoplus_{\sig \in \Sig(n)} \widehat{H}^*_\bT(x_\sig;\bQ)$$
obtained by restriction to $\bT$-fixed points 
is injective (e.g., see \cite[Subsection 2.5]{CMSSEM}). 
By localizing at the multiplicative set $L$, one gets by the exactness of localization an injective map
$$\bigoplus_{\sig \in \Sig(n)} i_\sig^* : \widehat{H}^*_\bT(X;\bQ)_L \hookrightarrow \bigoplus_{\sig \in \Sig(n)} \widehat{H}^*_\bT(x_\sig;\bQ)_L. $$
 So it is enough to check formula  \eqref{ieHirz0} by using,  for each fixed point $x_\sig$, $\sig \in \Sig(n)$, the induced restriction  map $pr_{x_\sig}$ 
 \be\label{indpr} \widehat{H}^*_\bT(X;\bQ) \to \widehat{H}^*_\bT(X;\bQ)_L  {\to}  \widehat{H}^*_\bT(x_\sig;\bQ)_L \to \widehat{H}^*_{\bT'}(x_\sig;\bC)_{L'},\ee
 with the middle arrow given by $\frac{i_\sig^*}{Eu^\bT_X(x_\sig)}$ as in \eqref{locmap}, and $L' \subset (\Lambda_{\bT'})_\bC=H^{*}_{\bT'}(pt;\bC)$ is the multiplicative set generated by the elements $\pm a \cdot c(m')$, for $0\neq m' \in M'$ and $a \in \bC^*$ (using the notations preceeding formula \eqref{isumsmusd}).
Also the direct sum $\bigoplus_{\sig \in \Sig(n)} pr_{x_\sig}$ of these induced restriction maps is still injective, since the localization map on the integral domain $\widehat{H}^*_\bT(x_\sig;\bQ) \to  \widehat{H}^*_\bT(x_\sig;\bQ)_L$ is injective, and $\frac{i_\sig^*}{Eu^\bT_X(x_\sig)}$
 differs from $i_\sig^*$ by the unit $$Eu^\bT_X(x_\sig)=
 \vert G_\sig \vert \prod_{\rho \in \sig(1)} i_\sig^*F_\rho \in \widehat{H}^*_\bT(x_\sig;\bQ)_L.$$
Moreover, no information is lost if we consider complex instead  of rational coefficients. 

It thus suffices to prove the equality of both sides of formula 
\eqref{ieHirz0} after applying the cohomological localization map  \eqref{locmap} for any $\sig \in \Sig(n)$, with $\frac{i_\sig^*T^\bT_{y*}(X)}{Eu^\bT_X(x_\sig)}= T^\bT_{y*}(X)_{x_\sig}$ as in   
\eqref{isum}. 

If $g \in G_\Sig \setminus G_\sig$, there exists a $\rho \in \Sigma(1)$ with $a_\rho(g) \neq 1$, so that one factor on the restriction of the right-hand side of \eqref{ieHirz0} to $U_\sig$ becomes $0$. Hence the summation on the right-hand side of \eqref{ieHirz0} reduces after restriction to $U_\sig$ to a summation over ${g \in G_\sig}$. If $g \in G_\sig$ and $\rho \notin \sig(1)$, then $a_\rho(g)=1$, so that the restriction of the factor $ \frac{ F_{\rho} \cdot 
\big( 1+y  \cdot a_{\rho}(g)  \cdot e^{-F_{\rho}}\big)}{1-a_{\rho}(g) \cdot e^{-F_{\rho}}} $ to $U_\sig$ becomes $1+y$. With $r-n=\vert \Sig(1) \setminus \sig(1)\vert$, this will in turn cancel the factor $(1+y)^{n-r}$ of formula \eqref{ieHirz0}. 
Therefore,  we get
\begin{align*}
&\frac{1}{Eu_X^\bT(x_\sig)} \cdot {i_\sig^*} \left( (1+y)^{n-r} \cdot \sum_{g \in G_{\Sig}}  
\prod_{\rho \in \Sigma(1)} \frac{ F_{\rho} \cdot 
\big( 1+y  \cdot a_{\rho}(g)  \cdot e^{-F_{\rho}}\big)}{1-a_{\rho}(g) \cdot e^{-F_{\rho}}}    \right)  \\
&=\frac{1}{Eu_X^\bT(x_\sig)} \cdot {i_\sig^*} \left(
\sum_{g \in G_{\sig}}  
\prod_{\rho \in \sig(1)} \frac{ F_{\rho} \cdot 
\big( 1+y  \cdot a_{\rho}(g)  \cdot e^{-F_{\rho}}\big)}{1-a_{\rho}(g) \cdot e^{-F_{\rho}}}    
\right) \\
&=\frac{1}{\vert G_\sig \vert} \sum_{g \in G_\sig}  {i_\sig^*}  \left(  
\prod_{\rho \in \sig(1) } \frac{ 
 1+y  \cdot a_{\rho}(g)  \cdot e^{-F_{\rho}} }{1-a_{\rho}(g) \cdot e^{-F_{\rho}}}    
 \right).
\end{align*}
This expression reduces to formula \eqref{isum}, as desired, after 
noticing that, for $\rho=\rho_i$ any ray of $\sig(1)$, one gets that $i_\sig^* F_{\rho_i}=c(m'_{\sig,i})$, as well as by changing $g$ by $g^{-1}$ in $G_\sig$.
\end{proof}

\subsection{Equivariant Hirzebruch classes of orbit closures}
Let $X=X_\Sig$ be as before a projective simplicial toric variety with torus $\bT=T_N$.
Then the fan of the toric variety $V_{\sigma}$ given by the closure of the orbit $O_\sig$ corresponding to the cone $\sigma \in \Sigma$ can be described as follows. Let $N(\sigma)=N/N_\sigma,$
with $N_\sigma$ denoting  the sublattice of $N$ spanned by the points in $\sigma \cap N$.  Let $T_{N(\sigma)}=N(\sigma) \otimes_\bZ \bC^*$ be the torus associated to  $N(\sigma)$. For each cone $\nu \in \Sigma$ containing $\sigma$, let $\overline{\nu}$ be the image cone in $N(\sigma)_\bR$ under the quotient map $N_\bR \to N(\sigma)_\bR$. Then 
\be \label{star}
Star(\sigma)=\{ \overline{\nu} \subseteq N(\sigma)_\bR \mid \sig \preceq \nu \}
\ee 
is a simplicial fan in $N(\sigma)_\bR$, with associated toric variety isomorphic to $V_{\sigma}$ (see, e.g., \cite[Prop.3.2.7]{CLS}). Note that $\bT$ acts on $V_{\sigma}$ via the morphism $\bT=T_N \to T_{N(\sigma)}$ induced by the quotient map $N \to N(\sigma)$.

The $\bT$-equivariant Hirzebruch class of the $\bT$-equivariant closed inclusion $i_\sig:V_\sig \hookrightarrow X$ of the projective simplicial toric variety $V_\sig$ is computed by  (see \cite[formula (238)]{CMSSEM}):
\be\label{dueqeeh}
T_{y*}^\bT([{V_\sigma} \hookrightarrow X])= (1+y)^{n - r} \cdot \sum_{g\in G_{Star(\sig)}} 
\mult(\sig) \cdot \prod_{\rho\in\sigma(1)}F_\rho \cdot
\prod_{\rho\in  Star(\sigma)(1)} \frac{ F_{\rho} \cdot 
\big( 1+y  \cdot a_{\rho}(g)  \cdot e^{-F_{\rho}}\big)}{1-a_{\rho}(g) \cdot e^{-F_{\rho}}}.
\ee
Here $\mult(\sig)=|G_\sig|$ is the multiplicity of $\sig$, and $\rho\in Star(\sigma)(1)$ is a short notation for $\rho \in \bigcup_{\sig \preceq \nu} \nu(1) \setminus \sig(1)$.
Also, $T_{y*}^\bT([{V_\sigma} \hookrightarrow X])=(i_\sig)_*T_{y*}^\bT(V_\sigma)$ by functoriality, and we have the identification $T_{y*}^\bT(V_\sigma)=T_{y*}^{T_{N(\sigma)}}(V_\sigma)$ via a choice of splitting for the projection $\bT \to T_{N(\sigma)}$ (see \cite[Corollary 3.13]{CMSSEM}). 
Formula \eqref{dueqeeh} follows from formula \eqref{ieHirz0} together with the projection formula and 
$$[V_\sig]_\bT=\mult(\sig) \cdot \prod_{\rho \in \sig(1)} F_\rho.$$
For $y=0$, one gets a similar formula for the pushforward of the equivariant Todd classes $$\td_*^\bT(V_\sigma)=T_{0*}^\bT(V_\sigma)=T_{0*}^{T_{N(\sigma)}}(V_\sigma)=\td_*^{T_{N(\sigma)}}(V_\sigma).$$

\section{Euler-Maclaurin formulae}
We now describe applications of the results mentioned above to Euler-Maclaurin formulae. 
\subsection{Euler-Maclaurin formulae via polytope dilatation}
Let $P$ be a full-dimensional lattice polytope in $M_\bR\simeq \bR^n$, with toric variety $X=X_P$, inner normal fan $\Sigma=\Sigma_P$, and ample Cartier divisor $D=D_P$. Let $\Sigma(1)$ be the set of rays of $\Sigma$, corresponding to the facets $F$ of $P$. For each ray $\rho \in \Sig(1)$, let $u_\rho \in N$ be the corresponding ray generator. As before, we let $F_\rho:=[D_\rho]_\bT$ be the equivariant fundamental class of the $\bT$-equivariant divisor $D_\rho$ on $X$ corresponding to the ray $\rho \in \Sig(1)$.
Let $P(h)$ be the \index{dilated polytope} dilation of $P$ with respect to the vector $h=(h_\rho)_{\rho \in \Sigma(1)}$ with real entries indexed by the rays of $\Sigma$. So, if $P$ is defined by inequalities of the form
$$\langle m, u_\rho\rangle +c_\rho \geq 0,$$
with $u_\rho$ the ray generators and $c_\rho \in \bZ$, for each $\rho \in \Sig(1)$, then $P(h)$  is defined by inequalities 
$$\langle m, u_\rho \rangle +c_\rho +h_\rho \geq 0,$$
for each $\rho \in \Sig(1)$. In these notations, we have that $D=D_P=\sum_{\rho \in \Sig(1)} c_\rho \cdot D_\rho$. If the $h_\rho$'s are small enough, then $P(h)$ is also a full-dimensional simple polytope.
Recall here that $\langle \cdot , \cdot \rangle: M\times N\to \bb{Z}$ is the canonical pairing. 

In the above notations, we have the following result from \cite[Theorem 5.1]{CMSSEM}.
\bt\label{them1} 
\be\label{f113}
\begin{split}
\int_{P(h)} e^{\langle  m, z\rangle} \ dm = \sum_{\sig \in \Sig(n)} \frac{ e^{\langle i_\sig^* [D_{P(h)}],z\rangle}}{\langle Eu^\bT_X(x_\sig),  z \rangle} 
=
 \sum_{\sig \in \Sig(n)} \frac{ e^{\langle i_\sig^* c^\bT_1(\cO_X(D_{P})),z\rangle}}{\langle Eu^\bT_X(x_\sig),  z \rangle} \cdot  e^{ \sum_\rho h_\rho \langle i_\sig^* F_\rho,z\rangle} \:,
 \end{split}
\ee
with the canonical pairing extended to $z \in N_\bC:=N \otimes_\bZ \bC=\Hom_\bR(M_\bR, \bC)$. 
\et

The proof of this theorem follows similar arguments of Brion-Vergne \cite[Thm.4.5]{BrV2}. 
Starting with the Riemann sum approximation of an integral, we have:
\be\label{rsum}\lim_{k\to \infty} \frac{1}{k^{n}}\sum_{m \in k\cdot P(h) \cap M} e^{\langle m, \frac{1}{k} \cdot z\rangle}= \lim_{k\to \infty} \frac{1}{k^{n}}\sum_{m \in P(h) \cap \frac{1}{k} M} e^{\langle  m, z\rangle} = \int_{P(h)} e^{\langle  m, z\rangle} \ dm,\ee
with the Lebesgue measure $dm$ normalized so that the unit cube in $M \subset M_\bR$ has volume $1$ (which explains the use of the factor $\frac{1}{k^n}$). 
If the $h_\rho$'s are small rational numbers, then one can choose a large $k \in \bN$ so that $k\cdot P(h)$ is a lattice polytope in $\bR^n$ with respect to the lattice $M$. The left hand side of \eqref{rsum} is then computed by applying the function $\langle -, \frac{1}{k} \cdot z \rangle$ to the localized Riemann-Roch formula (obtained by combining \eqref{ic1} $y=0$ with \eqref{f104}) for the lattice polytope $k\cdot P(h)$ with $k\cdot [D_{P(h)}]_\bT=[k \cdot D_{P(h)}]_\bT$ : 
\be\label{f110}
\sum_{m \in k\cdot P(h) \cap M} e^{s(m)}= \sum_{\sig \in \Sig(n)} \frac{i_\sig^* \left( e^{k\cdot [D_{P(h)}]_\bT}  \right)}{Eu^\bT_X(x_\sig)} \cdot i_\sig^*(\td^\bT_*(X))  \in  L^{-1}(\Lambda^{an}_\bT)_\bQ \subset L^{-1}(\widehat{\Lambda}_\bT)_\bQ \:,
\ee 
together with 
\be\label{to} \lim_{k\to \infty} \langle i^*_\sig (\td^\bT_*(X)), \frac{1}{k} \cdot z \rangle =1.\ee
The latter equality is derived by using the explicit formula \eqref{ieHirz0} for $\td^\bT_*(X)=T^\bT_{0*}(X)$.

\br\label{regular} The left hand side of \eqref{f113} is a {\it continuous} function in $h$  near zero, and for all $z \in N_\bC$, whereas the right hand side is an {\it analytic} function in $h$  near zero,
and for $z \in N_\bC$ away from the linear  hyperplanes  $\langle i_\sig^* F_\rho, z \rangle = 0$ for each ray $\rho\in \sigma(1)$ of $\sigma\in \Sig(n)$. But then both sides of this equality have to be 
{\it analytic  functions in $h$ near zero and all $z\in N_\bC$},  with the corresponding Taylor series around zero converging uniformly on small compact neighborhoods of zero in the variables $h$ and $z$ (cf. also \cite[page 27]{KSW}).
\er

Let $P$ be a full-dimensional simple lattice polytope in $M_\bR$, and fix a face $E$ of $P$. Let $\sigma:=\sigma_E$ be the corresponding cone in the inner normal fan $\Sigma=\Sigma_P$ of $P$, with $V_\sigma=V_{\sigma_E}=V_E$ the closure of the orbit of $\sigma$ in $X=X_P$. Denote by $i_E=i_{\sigma}: V_\sig \hookrightarrow X$ the closed inclusion map. Then 
$V_\sigma$ is a simplicial toric variety whose fan is $Star(\sigma)$, as defined in \eqref{star}, which is built from cones $\tau \in \Sigma$ that have $\sigma$ as a face.  With these notations we have similarly (see \cite[Theorem 7.12]{CMSSEM}):

\bt\label{them2} 
\be\label{f113b}
\int_{E(h)} e^{\langle  m, z\rangle} \ dm =
 {\rm mult}(\sig_E) \cdot \sum_{\sig \in \Sig(n)} \frac{ e^{\langle (i_\sig^* c^\bT_1(\cO_X(D_{P})),z\rangle}}{\langle Eu^\bT_X(x_\sig),  z \rangle} \cdot  e^{ \sum_\rho h_\rho \langle i_\sig^* F_\rho,z\rangle} \cdot  \prod_{\rho \in \sig_E(1)} \langle i_\sig^* F_\rho, z\rangle \:.
\ee
\et
Here, the Lebesgue measure $dm$ on $E(h)$ is normalized so that the unit cube in the lattice $Span(E_0)\cap M$ has volume $1$, with $E_0:=E-m_0$ a translation of $E$ by a vertex $m_0 \in E$.

\smallskip

Comparison of formulae \eqref{f113} and \eqref{f113b} one gets the following key formula relating integrals over the dilated polytope to integrals about the dilated faces:
\be\label{rel1}
\int_{E(h)} f(m) e^{\langle m, z \rangle}  \ dm = {\rm mult}(\sig_E) \cdot \prod_{\rho \in \sig_E(1)} \frac{\partial}{\partial h_\rho} \int_{P(h)} f(m) e^{\langle m, z \rangle}  \ dm \:,
\ee

Based on the above formulae, we get the following abstract Euler-Maclaurin formula coming from the equivariant Hirzebruch-Riemann-Roch theorem (see \cite[Theorem 5.18, Corollary 5.21. Proposition 5.22]{CMSSEM}):
\bt\label{abstrEMi}
Let $[\cF]\in K^\bT_0(X)$ be fixed, and choose a convergent power series $p(x_\rho) \in \bQ\{x_\rho \mid \rho \in \Sig(1) \}$ so that $p(F_\rho)=\td^\bT_*([\cF]) \in \left( H^*_\bT(X;\bQ) \right)^{an}$. Then, with $p(\frac{\partial}{\partial h})$ the corresponding infinite order differential operator obtained from $p(x_\rho)$ by substituting $x_\rho \mapsto \frac{\partial}{\partial h_\rho}$, for all $\rho \in \Sig(1)$, we have that, for any polynomial  function $f$ on $M_\bR$,
\begin{multline}\label{f1134i}
p(\frac{\partial}{\partial h}) \left( \int_{P(h)} f(m) \cdot e^{\langle  m, z\rangle} \ dm \right)_{\vert_{h=0}} =\\
= \sum_{m\in M} \left( \sum_{i=0}^n (-1)^i \cdot \dim_\bC H^i(X;\cO_X(D) \otimes \cF)_{\chi^{-m}}\right) \cdot f(m) \cdot e^{\langle m, z \rangle}
\:,
\end{multline}
as analytic functions in $z\in N_\bC$ with $z$ small enough.
\et
For $f=1$, one first proves the formula
\be\label{rot} 
\begin{split} p(\frac{\partial}{\partial h}) \left( \int_{P(h)} e^{\langle  m, z\rangle} \ dm \right)_{\vert_{h=0}}  &= \langle \chi^\bT(X,\cO_X(D) \otimes \cF) , z\rangle \\
&= \sum_{m\in M} \left( \sum_{i=0}^n (-1)^i \cdot \dim_\bC H^i(X;\cO_X(D) \otimes \cF)_{\chi^{-m}}\right) \cdot e^{\langle m, z \rangle}
\end{split}
\ee
via the application of the operator $p(\frac{\partial}{\partial h}) (-)\vert_{h=0}$ to formula \eqref{f113}, together with the localized equivariant Hirzebruch-Riemann-Roch formula \eqref{f104}. To get \eqref{f1134i}, we apply the operator $ f(\frac{\partial}{\partial z}) $ to the last term of formula \eqref{rot}, seen as a formal power series in $z$.

Note that by evaluating 
formula \eqref{f1134i} at $z=0$ and for $f=1$ (i.e., forgetting the $\bT$-action), we get a {\it generalized volume formula}, namely, 
\[
p(\frac{\partial}{\partial h}) \left( {vol} \ {P(h)}  \right)_{\vert_{h=0}} = \chi(X,\cO_X(D) \otimes \cF),
\]
with ${vol} \ {P(h)} =\int_{P(h)} dm$ the volume of $P(h)$ and the Lebesgue measure normalized so that the unit cube in $M \subset M_\bR$ has volume $1$. See \cite[Theorem 2.15]{BrV1} for the case when $\cF=\cO_X$ (corresponding to counting points in $P \cap M$) and $\cF=\omega_X$ (corresponding to counting points in $\Int(P) \cap M$).

\medskip

For suitable choices of $[\cF]\in K^\bT_0(X)$, formula \eqref{f1134i} can be specialized to yield old and new Euler-Maclaurin type formulae. We include below several such examples, but
see \cite[Sections 5.3 and 6.2]{CMSSEM} for more details and examples (like, e.g., 
a (weighted) Euler-Maclaurin formula for a simple lattice polytope with some facets removed, cf. \cite[formula (210)]{CMSSEM}).


\bex We list below several specializations of formula \eqref{f1134i}
for appropriate choices of $[\cF]\in K^\bT_0(X)$ and explicit convergent power series $p(x_\rho) \in \bQ\{x_\rho \mid \rho \in \Sig(1) \}$ so that $p(F_\rho)=\td^\bT_*([\cF]) \in \left( H^*_\bT(X;\bQ) \right)^{an}$. See \cite[Sections 5.3 and 6.2]{CMSSEM} for complete details.
\begin{enumerate}
\item[(a)] For $[\cF]=mC_y^\bT(X)\in K^\bT_0(X)[y]$, with corresponding operator given by formula \eqref{ieHirz0},
$${T}_{y}(\frac{\partial}{\partial h}):=(1+y)^{n-r} \cdot 
 \sum_{g \in G_{\Sig}}  \prod_{\rho \in \Sig(1)} \frac{ \frac{\partial}{\partial h_{\rho} } \cdot 
\big( 1+y  \cdot a_{\rho}(g)  \cdot e^{-\frac{\partial}{\partial h_{\rho} }}\big)}{1-a_{\rho}(g) \cdot e^{-\frac{\partial}{\partial h_{\rho} }}}, $$
formula \eqref{f1134i} becomes (for $z\in N_\bC$ small enough) the following weighted Euler-Maclaurin formula for $P$:
\be\label{wem0}
{T}_{y}(\frac{\partial}{\partial h}) \left( \int_{P(h)} f(m) \cdot e^{\langle m, z \rangle} dm \right)_{|_{h=0}} =\sum_{E \preceq P} (1+y)^{\dim(E)} \sum_{m \in \Relint(E) \cap M} f(m) \cdot e^{\langle m , z\rangle}.
\ee

\item[(b)] For $\cF=\cO_X$ the structure sheaf of $X$, with corresponding operator $Todd(\frac{\partial}{\partial h})={T}_{0}(\frac{\partial}{\partial h})$, formula \eqref{f1134i} reduces to the classical Euler-Maclaurin formula of Brion-Vergne \cite{BrV2} for simple lattice polytopes, and 
Khovanskii-Pukhlikov \cite{KP} for Delzant lattice polytopes 
(corresponding to smooth projective toric varieties), namely:
\be\label{itodd}
Todd(\frac{\partial}{\partial h}) \left( \int_{P(h)} f(m) \cdot e^{\langle  m, z\rangle} \ dm \right)_{|_{h=0}} = \sum_{m \in P \cap M} f(m) \cdot e^{\langle  m, z\rangle}.
\ee

\item[(c)] For $\cF=\omega_X$ the canonical sheaf of $X$, with corresponding {dual Todd operator}
\be
Todd^\vee(\frac{\partial}{\partial h}):= \sum_{g\in G_{\Sigma}} \prod_{\rho\in \Sigma(1)} 
\frac{ a_{\rho}(g) \cdot \frac{\partial}{\partial h_{\rho}}\cdot e^{-\frac{\partial}{\partial h_{\rho}} }}{1-a_{\rho}(g)\cdot e^{-\frac{\partial}{\partial h_{\rho}} } }
\in \bb{Q}\{\frac{\partial}{\partial h_{\rho} }\;|\; \rho \in \Sigma(1)\} \:,
\ee
given by the coefficient of $y^n$ in ${T}_{y}(\frac{\partial}{\partial h})$, formula \eqref{f1134i} reduces to the Euler-Maclaurin formula for the interior of $P$:
\be\label{tdu}
Todd^\vee(\frac{\partial}{\partial h}) \left( \int_{P(h)} f(m) \cdot e^{\langle  m, z\rangle} \ dm \right)_{|_{h=0}} 
 = \sum_{m\in \Int(P)\cap M} f(m) \cdot e^{\langle  m, z\rangle} \:.
\ee
\item[(d)] Let $[\cF]=[(i_E)_*mC_y^\bT(V_E)]\in K^\bT_0(X)[y]$, with $E$ a face of $P$ and $i_E=i_\sigma: V_{\sigma}=V_E \hookrightarrow X$ the inclusion of the orbit closure for the cone $\sigma$ corresponding to $E$. 
For the operator given by formula \eqref{dueqeeh}, i.e., 
\be
T_y^E(\frac{\partial}{\partial h}):=(1+y)^{n - r} \cdot \sum_{g\in G_{Star(\sig)}} \mult(\sig) \cdot \prod_{\rho\in\sigma(1)}\frac{\partial}{\partial h_\rho} \cdot \prod_{\rho\in  Star(\sigma)(1)} 
\frac{  \frac{\partial}{\partial h_{\rho}} \left( 1+y\cdot a_{\rho}(g)\cdot e^{-\frac{\partial}{\partial h_{\rho}} } \right)}
{1-a_{\rho}(g)\cdot e^{-\frac{\partial}{\partial h_{\rho}} } }.
\ee
formula \eqref{f1134i} reduces to a weighted Euler-Maclaurin formula for the face $E$ of $P$:
\be\label{EMinth}
T_y^E(\frac{\partial}{\partial h})\left( \int_{P(h)} f(m) \cdot e^{\langle  m, z\rangle} \ dm \right)_{|_{h=0}} =
\sum_{E' \preceq E} (1+y)^{\dim(E')} \cdot \sum_{m \in \Relint(E') \cap M} f(m)\cdot e^{\langle  m, z\rangle}.
\ee
\item[(e)] For $[\cF]=[(i_E)_*\cO_{V_E}]\in K^\bT_0(X)$, with operator 
$Todd_E(\frac{\partial}{\partial h})={T}^E_{0}(\frac{\partial}{\partial h})$,   
formula \eqref{f1134i} reduces to a Euler-Maclaurin formula for the face $E$ of $P$:
\be\label{EMint}
Todd_E(\frac{\partial}{\partial h})\left( \int_{P(h)} f(m)\cdot e^{\langle  m, z\rangle} \ dm \right)_{|_{h=0}} =\sum_{m \in E \cap M} f(m)\cdot e^{\langle  m, z\rangle},
\ee
\item[(f)] For $[\cF]=[(i_E)_*\omega_{V_E}]\in K^\bT_0(X)$, with operator
\be
Todd^\vee_E(\frac{\partial}{\partial h}):= \sum_{g\in G_{Star(\sig)}} \mult(\sig) \cdot \prod_{\rho\in\sigma(1)}\frac{\partial}{\partial h_\rho} \cdot \prod_{\rho\in  Star(\sigma)(1)} \frac{ a_{\rho}(g) \cdot \frac{\partial}{\partial h_{\rho}}\cdot e^{-\frac{\partial}{\partial h_{\rho}} }}{1-a_{\rho}(g)\cdot e^{-\frac{\partial}{\partial h_{\rho}} } }
\ee
given by the coefficient of $y^{\dim(E)}$ in ${T}^E_{y}(\frac{\partial}{\partial h})$,
formula \eqref{f1134i} reduces to a Euler-Maclaurin formula for the interior of the face $E$ of $P$:
\be\label{EMinta}
Todd^\vee_E(\frac{\partial}{\partial h})\left( \int_{P(h)} f(m)\cdot e^{\langle  m, z\rangle} \ dm \right)_{|_{h=0}} =\sum_{m \in \Relint(E) \cap M} f(m)\cdot e^{\langle  m, z\rangle},
\ee
\end{enumerate}
\qed
\eex

\bex
More general examples of explicit weighted Euler-Maclaurin formulae can be obtained by 
twisting the classes $[\cF]:=[mC_y^\bT(X)] \in K_0^\bT(X)[y]$
by $\cO_X(D'-D)$, for $D=D_P$ the original ample divisor associated to the full-dimensional simple lattice polytope $P$, and $D'$ any $\bT$-invariant Cartier divisor on $X$ (see \cite[Theorem 6.2, Corollary 6.3, Example 6.9]{CMSSEM}).

Let $D'$ be a globally generated $\bT$-invariant Cartier divisor on $X$, with associated (not necessarily full-dimensional) lattice polytope $P_{D'} \subset M_\bR$, e.g., as in \cite[Prop.6.2.3]{CLS}. Let $D'-D=\sum_{\rho \in \Sig(1)} d_\rho D_\rho$ as a $\bT$-invariant Cartier divisor. 
Let $X_{D'}$ be the toric variety of the lattice polytope $P_{D'}$, defined via the corresponding {generalized fan}. 
Consider the infinite order differential operator
\begin{equation}\label{itodd222b}
T'_y(\frac{\partial}{\partial h}):= e^{\sum_{\rho \in \Sig(1)} d_\rho \cdot \frac{\partial}{\partial h_{\rho} }} \cdot T_y(\frac{\partial}{\partial h})  \in \bQ\{ \frac{\partial}{\partial h_\rho} \mid \rho \in \Sig(1)\}[y].
\end{equation}
Formula \eqref{f1134i} reduces to the following new weighted Euler-Maclaurin formula:
\begin{multline}\label{iwem01b}
T'_{y}(\frac{\partial}{\partial h}) \left( \int_{P(h)} f(m) \cdot e^{\langle  m, z\rangle} \ dm \right)_{|_{h=0}} = \\ =\sum_{E \preceq P_{D'}} \left( \sum_{\ell \geq 0} 
(-1)^\ell \cdot d_\ell(X/E) \cdot (1+y)^{\ell + \dim(E)}   \right) \cdot \sum_{m \in \Relint(E) \cap M} f(m) \cdot e^{\langle  m, z\rangle},
\end{multline}
with multiplicities $d_\ell(X/E)=d_\ell(X/\sig')$ as in \eqref{bl1}, and the face $E$ of $P_{D'}$ corresponding to the cone $\sig' \in \Sig'$.
Note that in this context $P_{D'}$ is an $\bN$-Minkowski summand of the original polytope $P$, see \cite[Cor.6.2.15]{CLS}.
Forgetting the $\bT$-action (i.e., for $f=1$ and $z=0$), one gets the following volume formula (fitting with \eqref{fort}):
\be
T'_{y}(\frac{\partial}{\partial h}) \left( vol \ P(h) \right)_{\vert_{h=0}}=
\sum_{E \preceq P_{D'}} \left( \sum_{\ell \geq 0} 
(-1)^\ell \cdot d_\ell(X/E) \cdot (1+y)^{\ell + \dim(E)}  \right) \cdot  \vert \Relint(E) \cap M \vert .
\ee
\qed
\eex

\subsection{Euler-Maclaurin formulae of Cappell-Shaneson type}
Let $\sig_E$ be the cone in $\Sig=\Sig_P$ corresponding to the face $E$ of $P$, and let $V_E=V_{\sig_E}$ be the closure of the $\bT$-orbit in $X$ corresponding to $\sig_E$. Then the equivariant fundamental classes $[V_{E}]_\bT$ generate $\widehat{H}^*_\bT(X;\bQ)$ as a  
$\widehat{H}^*_\bT(pt;\bQ)= (\widehat{\Lambda}_\bT)_\bQ$-algebra (e.g., see \cite[Remark 7.15]{CMSSEM}).
Let now $[\cF]\in K^\bT_0(X)$ be fixed, and choose 
elements $p_{E}(t_i)\in \widehat{H}^*_\bT(pt;\bQ)= (\widehat{\Lambda}_\bT)_\bQ \overset{s}{\simeq} \bQ[[t_1,\dots,t_n]]$ with
\begin{equation}\label{main-toddi}
\td_*^\bT([\cF])= \sum_{E\preceq P} \: p_{E}(t_i)[V_{E}]_\bT \in \widehat{H}^*_\bT(X;\bQ)\:.
\end{equation}
With these notations, we can now state an abstract Euler-Maclaurin formula coming from the equivariant Hirzebruch-Riemann-Roch theorem, in terms of integrals over the faces of a polytope, instead of using a dilated polytope, see \cite[Theorem 7.16]{CMSSEM}:
\bt\label{abstrEM3i}
Let $X=X_P$ be the projective simplicial toric variety associated to the full-dimensional simple lattice polytope $P\subset M_\bR$. Let $\Sig:=\Sig_P$ be the inner normal fan of $P$, and $D:=D_P$ the ample Cartier divisor associated to $P$. 
Then for a polynomial  function $f$ on $M_\bR$, we have:
\be\label{cs2137i}
\sum_{E\preceq P} \:  \int_{E} p_{E}(\partial_i) f(m) \ dm 
= \sum_{m\in M} \left( \sum_{i=0}^n (-1)^i \cdot \dim_\bC H^i(X;\cO_X(D) \otimes \cF)_{\chi^{-m}}\right) \cdot f(m)\:.
\ee
\et
Here, $\partial_i=\frac{\partial}{\partial t_i}$, with respect to the coordinates $t_i$ of $M_\bR\simeq \bR^n$.
Theorem \ref{abstrEM3i} follows by applying the operator $p(\partial_i, \frac{\partial}{\partial h})(-)|_{h=0}$ with 
\be\label{main-todd2}
p(\partial_i, \frac{\partial}{\partial h}):=\sum_{E\preceq P} \: \left( {\rm mult}(\sig_E) \cdot \prod_{\rho \in \sig_E(1)} \frac{\partial}{\partial h_\rho}\right) \cdot p_{\sig_E}(\partial_i), 
\ee
 to formula \eqref{f113}, together with the key formula \eqref{rel1} and the localized equivariant Hirzebruch-Riemann-Roch formula \eqref{f104}.


In the classical case $\cF:=\cO_X$, 
this is exactly Cappell-Shaneson's recipe for the definition of the differential operators $p_{E}(\partial_i)$, described here geometrically in terms of 
the equivariant Todd class $\td_*^\bT(X):=\td_*^\bT([\cO_X])\in  \widehat{H}^*_\bT(X;\bb{Q})$  (see \cite[Theorem 2]{CS2} or \cite[Section 6.2]{S}). 
In this case, \eqref{cs2137i} reduces to the Cappell-Shaneson Euler-Maclaurin formula:
\be
\sum_{E\preceq P} \:  \int_{E} p_{E}(\partial_i) f(m) \ dm 
= \sum_{m\in P\cap M} f(m).
\ee
See \cite[Example 7.18]{CMSSEM} for further specializations of formula \eqref{cs2137i}.

\subsubsection{Generalized reciprocity for Dedekind sums via Euler-Maclaurin formulae}
We conclude this note with the following application of formula \eqref{cs2137i}, see \cite[Corollary 7.19]{CMSSEM}: 

\bc\label{coro8} In the context of Theorem \ref{abstrEM3i}, 
one gets the following identity:
\be\label{274}
\sum_{v \in P} \left( p_{{v}}(\partial_i) f \right) (0) = \sum_{m\in M} \left( \sum_{i=0}^n (-1)^i \cdot \dim_\bC H^i(X;\cF)_{\chi^{-m}}\right) \cdot f(m).
\ee
where the left hand sum is over the vertices of $P$.
\ec

\bex
If $\cF=\cO_X$ in \eqref{274}, one gets the following identity:
\be
\sum_{v \in P} \left( p_{{v}}(\partial_i) f \right) (0) = f(0).
\ee
For instance, in the case of lattice polygones, this formula yields generalizations of \index{reciprocity law} reciprocity laws for classical \index{Dedekind sum} {\it Dedekind sums} (using, e.g., the explicit description of the operators $p_{{v}}(\partial_i)$ from \cite[page 889]{CS2}). 
\qed
\eex

\begin{ack}
L. Maxim 
acknowledges support from the project ``Singularities and Applications'' - CF 132/31.07.2023 funded by the European Union - NextGenerationEU - through Romania's National Recovery and Resilience Plan.
J. Sch\"urmann was funded by the Deutsche Forschungsgemeinschaft (DFG, German Research Foundation) Project-ID 427320536 -- SFB 1442, as well as under Germany's Excellence Strategy EXC 2044 390685587, Mathematics M\"unster: Dynamics -- Geometry -- Structure. L. Maxim and  J. Sch\"urmann also thank the Isaac Newton Institute for Mathematical Sciences for the support and hospitality during the program ``Equivariant methods in geometry'' when work on this paper was undertaken.
\end{ack}





\printindex

\end{document}